\documentclass[10pt]{article}
\usepackage[english]{babel}
\usepackage{amsfonts}
\usepackage[dvips]{graphicx}

\begin{document}

\title{\bf Costruction of classic exact solutions for Tricomi equation}

\date{2010, October}

\author{\bf Gianluca Argentini \\
\normalsize{Research \& Development Dept., Riello Burners - Italy}\\
\normalsize gianluca.argentini@rielloburners.com \\
\normalsize gianluca.argentini@gmail.com \\}

\maketitle

\noindent{\bf Abstract}\\
A formula to construct classic exact solutions to Tricomi partial differential equation. The steps to obtain this formula require only elementary resolution of a simple system of first order PDEs.\\


The (generalized) Tricomi equation is the second order partial differential equation

\begin{equation}\label{tricomiPDE}
	\partial_{xx}f + a(x)\partial_{yy}f = 0
\end{equation}

\noindent where $f(x,y)$ and $a(x)$ are regular functions. In the case $a(x)=x$, (\ref{tricomiPDE}) is the classic Tricomi equation, as described by Prof. Francesco in his work dated 1923 \cite{tricomi}. The equation is an abstraction of the Euler equation on a 2D fluid motion in the case of a flow speed near the sonic condition. For Tricomi equation some exact solutions formulas exist (\cite{polyanin} or \cite{eqworld}), and solutions to some particular boundary values problems are known \cite{cibrario}. Also, some results about weak solutions are found \cite{yagdjian}.

Let be $u = \partial_x f = f_x$, $v = \partial_y f = f_y$. Then Tricomi equation can be written as a system of two first order equations:

\begin{equation}
\left\{
\begin{array}{ll}\label{systemTricomi}
	u_x(x,y) + a(x)v_y(x,y) = 0\\
	u_y(x,y) - v_x(x,y) = 0
\end{array}
\right.
\end{equation}

\noindent We define, for convenience reasons, a function $t = t(x,y) = v(x,y)$. From the first equation we can write

\begin{equation}\label{u_first}
	u(x,y) = - \int_a^x \left[ a(s) \hspace{0.1cm} t_y(s,y) \right] ds + g(y)
\end{equation}

\noindent where $a \in \mathbb{R}$ and $g$ is an arbitrary function of real variable. Now we have to find the condition that $g$ and $t$ must satisfy to verify the second equation of the system (\ref{systemTricomi}). This equation is now, from the usual rule of derivation of integrals depending on parameters,

\begin{equation}
	- \int_a^x \left[ a(s) \hspace{0.1cm} t_{yy}(s,y) \right] ds + g'(y) = t_x
\end{equation}

\noindent and hence

\begin{equation}\label{g_first}
	g(y) = \int_b^y \left[ t_x(x,r) + \int_a^x \left[ a(s) \hspace{0.1cm} t_{yy}(s,r) \right] ds \right] dr
\end{equation}

\noindent with $b \in \mathbb{R}$. But $g$ must depend only on $y$: using the Fundamental Theorem of the Integral Calculus we have

\begin{equation}
	0 = \partial_x g(y) = \int_b^y \left[ t_{xx}(x,r) +  a(x) \hspace{0.1cm} t_{yy}(x,r) \right] dr \hspace{0.5cm} \forall y
\end{equation}

\noindent If $t$ is a solution of Tricomi equation (\ref{tricomiPDE}), previous condition is verified. Note that, integrating with respect to $y$ in $g(y)$ formula, the expression for $u$ becomes

\begin{equation}\label{u_second}
	u(x,y) = \int_b^y t_x(x,r) dr - \int_a^x \left[a(s) \hspace{0.1cm} t_y(s,b) \right] ds
\end{equation}

\noindent So, we have found that if $t$ is a solution of Tricomi equation, then a function $f$ such that $f_x = u$ and $f_y = t$, with $u$ given by (\ref{u_second}), is a solution too. Hence, from $f_y = t$ we can write

\begin{equation}
	f(x,y) = \int_b^y t(x,r) dr + h(x)
\end{equation}

\noindent where $h$ is an arbitrary function of real variable. From $f_x = u$ follows

\begin{equation}
	\int_b^y t_x(x,r) dr + h'(x) = \int_b^y t_x(x,r) dr - \int_a^x \left[a(s) \hspace{0.1cm} t_y(s,b) \right] ds
\end{equation}

\noindent so that 

\begin{equation}
	h(x) = - \int_a^x \int_a^q \left[a(s) \hspace{0.1cm} t_y(s,b) \right] ds \hspace{0.1cm} dq
\end{equation}

\noindent We have so proven that {\it if $t=t(x,y)$ is a generic solution of Tricomi equation, then the following formula gives a solution of the same equation too}:

\begin{equation}\label{formulaSolution}
	f(x,y) = \int_b^y t(x,r) dr - \int_a^x \int_a^q \left[a(s) \hspace{0.1cm} t_y(s,b) \right] ds \hspace{0.1cm} dq
\end{equation}\\

\noindent For {\it a-posteriori} verification, let be $t$ such that $t_{xx}+a(x)t_{yy}=0$. We have

\begin{eqnarray}
	f_{xx} = \int_b^y t_{xx}(x,r) dr - a(x) \hspace{0.1cm} t_y(x,b) = \nonumber \\
	= - \int_b^y a(x)t_{yy}(x,r) dr - a(x) \hspace{0.1cm} t_y(x,b) = \nonumber \\
	= - a(x) \hspace{0.1cm} t_y(x,y) + a(x) \hspace{0.1cm} t_y(x,b) - a(x) \hspace{0.1cm} t_y(x,b) = \nonumber \\
	= - a(x) \hspace{0.1cm} t_y(x,y) = \nonumber \\
	= -a(x)f_{yy}
\end{eqnarray}	

\noindent so the solving formula (\ref{formulaSolution}) is verified.\\

{\it Example} 1. Let be $f_{xx} + x f_{yy} = 0$ the original Tricomi equation. Note that $t(x,y) = y$ is a (trivial) solution. Then, from (\ref{formulaSolution}) with $a=b=0$, the function

\begin{equation}
	f(x,y) = \int_0^y y \hspace{0.1cm} dr - \int_0^x \int_0^q s \hspace{0.1cm} ds \hspace{0.1cm} dq = \frac{1}{2}y^2 - \frac{1}{6}x^3
\end{equation}

\noindent is a solution too.\\

{\it Example} 2. Let be $f_{xx} + \textnormal{cos}(x) f_{yy} = 0$ a generalized Tricomi equation. Note that $t(x,y) = y$ is a (trivial) solution. Then, from (\ref{formulaSolution}) with $a=b=0$, the function

\begin{equation}
	f(x,y) = \int_0^y y \hspace{0.1cm} dr - \int_0^x \int_0^q \textnormal{cos}(s) \hspace{0.1cm} ds \hspace{0.1cm} dq = -1 + \frac{1}{2}y^2 + \textnormal{cos}(x)
\end{equation}

\noindent is a solution too.\\

\noindent Note that formula (\ref{formulaSolution}) can be used as a {\it solutions machine}: starting from the most simple not null solution, that is from $t(x,y) = 1$, it can be used for iterative construction of solutions.\\

Also, formula (\ref{formulaSolution}) can be used to find solutions to some particular boundary values problems for Tricomi equation. For example, in the case $a(x)=x$ and $b=0$, if the function $t(x,y)$ is such that $t_y(x,0)=0$ $\forall x$ (condition satisfied e.g. by $t=-\frac{1}{6}x^3+\frac{1}{2}y^2$), then from (\ref{formulaSolution}) follows

\begin{equation}
	f(x,y) = \int_0^y t(x,r) dr
\end{equation}

\noindent so that $f(x,0)=0$ $\forall x\in \mathbb{R}$. Hence, if $t_y(x,0)=0$ $\forall x$, then $f(x,y) = \int_0^y t(x,r) dr$ {\it is a solution of the boundary values differential problem}

\begin{equation}
	\partial_{xx}f + a(x)\partial_{yy}f = 0 \hspace{0.1cm}, \hspace{0.5cm} f(x,0) = 0	\hspace{0.3cm} x\in \mathbb{R}
\end{equation}

\noindent If $t$ is a solution of (\ref{tricomiPDE}) satisfying the boundary conditions $t(x,0)=g(x)$, with $g$ a function of real variable, then from (\ref{formulaSolution}) follows that $f_y(x,0)=g(x)$, that is {\it $f$ is a solution of the problem with Cauchy-Neumann condition}

\begin{equation}
	\partial_{xx}f + a(x)\partial_{yy}f = 0 \hspace{0.1cm}, \hspace{0.5cm} \partial_y f(x,0) = g(x)	\hspace{0.3cm} x\in \mathbb{R}
\end{equation}

\end{document}